\documentstyle[12pt]{article}
\newcommand{\inverse}[1]{{\textstyle\frac{1}{#1}}}
\newcommand{\half}{\inverse{2}}

\newcommand{\qed}{\hfill$\Box$}

\title{
{\Large \bf A finite dimensional filter\\ with exponential conditional density
}\thanks{Available also at www.damianobrigo.it. Part of this paper has been published in "Statistics and Probability Letters" 49 (2000), pp. 127-134.}
}
\author{Damiano Brigo\thanks{A short version co-athored by Francois Le Gland has appeared on the proceedings of the conference 36th IEEE CDC San Diego,  {\tt http://www.irisa.fr/sigma2/legland/pub/cdc97-fdf.pdf}}  \\
Department of Risk Management\\
CARIPLO Bank\\
via Boito 7\\
20121 Milano, Italy\\
e-mail : dbrigo@opoipi.it\\
}
\date{July 1997}
\newtheorem{theorem}{Theorem}[section]

\newtheorem{definition}[theorem]{Definition}
\newtheorem{example}[theorem]{Example}
\newtheorem{remark}[theorem]{Remark}
\newtheorem{problem}[theorem]{Problem}
\begin{document}
%
%\PARstart{X}{YYY} ZZZ

\maketitle
\thispagestyle{empty}
\begin{abstract}
In this paper we consider the continuous--time nonlinear filtering problem,
which  has an infinite--dimensional solution in general, as proved 
by Chaleyat--Maurel and Michel. 
There are few examples of nonlinear systems for which the optimal filter
is finite dimensional, in particular Kalman's, Bene\v{s}', and
Daum's filters. 
In the present paper, we construct new classes of scalar
nonlinear filtering problems 
admitting finite--dimensional filters.  
We consider a given (nonlinear) diffusion coefficient
for the state 
equation, a given (nonlinear) observation function, and a given 
finite--dimensional exponential family of probability densities. 
We construct a drift for the state equation such that the resulting
nonlinear filtering problem admits a finite--dimensional filter 
evolving in 
the prescribed exponential family augmented by the observaton function
and its square.  
\end{abstract}

{\bf keywords} Scalar Nonlinear Diffusion Processes,
Finite Dimensional Families,
Exponential Families,
Stochastic Differential Equations, 
Scalar Nonlinear Filtering Problem, Finite--Dimensional Filters.

\section{Introduction}
%This paper moves from the differential geometric approach to
%nonlinear filtering as developed by Brigo, Hanzon and LeGland
%(1995 --96) \cite{BrHaLe}, \cite{BrHaLe2} and \cite{Brigo2}. 
In this paper we consider the (scalar) nonlinear filtering problem 
in continuous time. 
For a quick introduction to the filtering problem see
Davis and Marcus (1981) \cite{davis81b}. For a more complete
treatment see Liptser and Shiryayev (1978) \cite{LipShi}
from a mathematical point of view or Jazwinski (1970) \cite{jazwinski70a}
for a more applied perspective.  

The nonlinear filtering problem has an infinite--dimensional solution
in general. Constructing nonlinear systems for which the optimal filter
is finite dimensional is a problem which received considerable 
attention in the past.
It turned out that such systems are quite rare. Examples were given
by Bene\v{s} \cite{benes81a} and Daum \cite{daum86a}. Instead, general results
on nonexistence of such systems, based on Lie--algebraic techniques,
were made available by Chaleyat--Maurel and Michel (1984) \cite{ChalMaur},
and related works appeared for example in 
Ocone and Pardoux (1989) \cite{OconPard}, 
L\'evine (1991) \cite{Levi}.
In the present paper, we construct scalar nonlinear filtering problems 
admitting finite--dimensional filters. Our examples can be added to the 
contributions of Bene\v{s} and Daum of known examples of finite dimensional
filters. 
The solution of the filtering problem is a Stochastic PDE
which can be seen as a generalization of the Fokker--Planck
equation expressing the evolution of the 
density of a diffusion process. 
This filtering equation is called Kushner--Stratonovich equation,
and an unnormalized (simpler) version of it is known as
the Duncan--Mortensen--Zakai Stochastic Partial Differential
Equation (DMZ--SPDE). It is this second equation that we shall consider
in the present article.
In \cite{BrHaLe}
%\cite{BrHaLe2} 
and \cite{BrigoPhD}
the Fisher metric is used to
project the Kushner--Stratonovich (or the Fokker--Planck) equation 
onto an exponential family of probability densities,
yielding the new class of approximate filters called {\em projection
filters}. In \cite{Brigo2,Brigo3} the Gaussian projection filter
is studied in the small-noise setting.
The results given in the present paper originate from ideas
expressed in such past works, especially in \cite{BrigoPhD} and
\cite{BrigoCNR}.
 
Our approach to the construction of nonlinear filtering
problems admitting finite--dimensional filters is the following. 
We model the state according to a SDE whose drift may depend on the 
observations. We consider a given (nonlinear) diffusion coefficient
for the state 
equation, a given (nonlinear) observation function, and a given 
finite--dimensional exponential family of probability densities. 
We construct a drift for the state equation such that the resulting
 nonlinear filtering problem admits a finite--dimensional filter 
evolving in 
the prescribed exponential family augmented in a particular way 
(depending on $h$).  

%Hazewinkel, Marcus and Sussmann (1983) \cite{HaMaSu},
 
\section{Problem formulation}
We start by introducing the filtering problem for 
scalar continuous time systems.

On the probability space $(\Omega,{\cal F},P)$ 
with the filtration $\{{\cal F}_{t}\,,\,t \in [0 \ T] \ \}$ 
we consider the following scalar state and observation equations:
%%  The case where the drift f_t depends on the whole history of the
%%  observations {\cal Y}_t up to time t and not only on Y_t:
%%
%%\begin{equation} \label{Lanc1-1} 
%%\begin{array}{rcl} 
%%   dX_t &=& f_t(X_t,Z_t)\,dt + \sigma_t(X_t)\,dW_t, \ \ X_0, \\ \\ 
%%   dY_t &=& h_t(X_t)\,dt + dV_t, \ \ Y_0 = 0\ . 
%%\end{array} 
%%\end{equation} 
%%The process $Z$ is adapted to the filtration 
%%${\cal Y}_t:=\sigma(Y_s\,,\,0\leq s\leq t)$.
%%
\begin{equation} \label{Lanc1-1} 
\begin{array}{rcl} 
   dX_t &=& f_t(X_t,Y_t)\,dt + \sigma_t(X_t)\,dW_t, \ \ X_0, \\ \\ 
   dY_t &=& h(X_t)\,dt + dV_t, \ \ Y_0 = 0\ . 
\end{array} 
\end{equation} 
We set 
\begin{displaymath}
a_t(\cdot) := \sigma_t(\cdot)^2 \ .
\end{displaymath}
Time invariance of $h$ is needed to simplify exposition. The reason why
we require it will result clear later on.
These equations are It\^o stochastic differential equations (SDE's). 
We shall use both It\^o SDE's (for example for the
signal $X$) and McShane--Fisk--Stratonovich (MFS) SDE's 
(when dealing with densities). The MFS form
will be denoted by the presence of the symbol `$\circ$' in between
the diffusion coefficient and the Brownian motion of a SDE. 
The noise processes $\{W_t\,,\,t \in [0 \ T] \ \}$ and $\{V_t\,,\,t \in [0 \ T] \ \}$ 
are two standard Brownian motions. 
Finally, the initial state $X_0$ and the noise processes 
$\{W_t\,,\,t \in [0 \ T] \ \}$ and $\{V_t\,,\,t \in [0 \ T] \ \}$ are assumed to be
independent. 

Notice that this model is different from the models given usually
(Jazwinski~\cite{jazwinski70a}, %Maybeck~\cite{maybeck79b}, 
Davis and Marcus~\cite{davis81b}) due to the presence of  
%%   the ${\cal Y}_t$--measurable term $Z_t$
$Y_t$ in the drift $f_t$ of the state equation.
 However, this does not complicate matters. 
Indeed, in \cite{LipShi} a general formulation is given of which
our model is a particular case. 
%%Moreover, in stochastic control such models are well known:
%%Consider for example the model (8.1.5), (8.1.10) in the book
%%\cite{bensoussan92b}. 
See also, for example, the application of SPDEs theory to 
filtering \cite{gyongy90a}, where
the used model class allows dependence of the coefficients $f,a$ and $h$
on the observation process $Y$. Such models are often encountered in
stochastic control, see for example the model (8.1.5), (8.1.10) in 
\cite{bensoussan92b}. 
The nonlinear filtering problem consists in finding 
the conditional probability distribution $\pi_t$ of the state $X_t$ 
given the observations up to time $t$, 
i.e.\ $\pi_t(dx) := P[X_t\in dx\mid {\cal Y}_t]$, 
where ${\cal Y}_t:=\sigma(Y_s\,,\,0\leq s\leq t)$. 
We shall say that a filtering problem such as (\ref{Lanc1-1})
satisfies condition $(A)$ if:
\begin{itemize}
\item[(A)]
 for all $t \in [0 \ T]$, 
the probability distribution $\pi_t$ has an unnormalized density $q_t$ 
w.r.t.\ the Lebesgue measure,
and $\{q_t\,,\,t \in [0 \ T] \ \}$ satisfies the It\^o--type stochastic
partial differential equation (Duncan--Mortensen--Zakai SPDE) 
\begin{equation} \label{DMZ} 
   dq_t = {\cal L}_t^\ast\, q_t\,dt 
   + q_t\, h \  dY_t  
\end{equation} 
in a suitable functional space, 
%where $E_{p_t}\{\cdot\}$ denotes the expectation w.r.t.\ 
%the probability density $p_t$, i.e.\ the conditional expectation 
%given the observations up to time $t$, 
%and 
where for all $t \geq 0$, 
the forward diffusion operator ${\cal L}_t^\ast$ is defined by 
\begin{displaymath} 
   {\cal L}_t^\ast \phi = - 
   \frac{\partial}{\partial x}\, [f_t\, \phi] 
   + \half 
   \frac{\partial^2}{\partial x^2}\, 
   [a_t\, \phi]\ , 
\end{displaymath} 
for any test function $\phi$.
\end{itemize}
Conditions under which $(A)$ holds, and more general results on
solutions of the DMZ--SPDE can be found in \cite{gyongy90a}. 
The corresponding MFS form of the SPDE~(\ref{DMZ}) is~:  
\begin{displaymath} 
   dq_t = {\cal L}_t^\ast\, q_t\,dt 
   - \half\, q_t\, \vert h \vert^2   \,dt 
   +   q_t\, h \ \circ dY_t\ . 
\end{displaymath} 
%In order to simplify notation, we introduce the following 
%definition, which will be used throughout the following sections~: 
%\begin{eqnarray} \label{coeff}
%   && \alpha_t(p) :=  
%   \frac{{\cal L}_t^\ast\, p}{p}\ =
%   - \sum_{i=1}^n [\, 
%   f_t^i\, \frac{\partial}{\partial x_i}(\log p) 
%   + \frac{\partial f_t^i}{\partial x_i} \,] 
%   \nonumber \\ \label{coeff:alpha} \\ \nonumber 
%   && \mbox{} + \half \sum_{i,j=1}^n [\, 
%   a_t^{ij}\, 
%   \frac{\partial^2}{\partial x_i\partial x_j}(\log p)  
%   + a_t^{ij}\, \frac{\partial}{\partial x_i}(\log p)\, 
%   \frac{\partial}{\partial x_j}(\log p) \\ \nonumber \\ \nonumber
%   && \hspace{2.5cm}  + 2\, \frac{\partial a_t^{ij}}{\partial x_j}\, 
%   \frac{\partial}{\partial x_i}(\log p) 
%   + \frac{\partial^2 a_t^{ij}}{\partial x_i\partial x_j} \,]\ . 
%\end{eqnarray} 
%
It is known that in general the density $q_t$ does not evolve in a
finite--dimensional parametrized family, say 
$\{q(\cdot,\zeta), \ \ \zeta \in U \subset {\bf R}^{m+1}\}$.   
In some particular cases this happens. For example, in the linear
case ($f$ linear, $\sigma_t(x) = \sigma_t$ for all $x$, 
$h$ linear, and $X_0$ with Gaussian distribution) 
$q_t$ evolves in the manifold of the unnormalized
Gaussian densities. Some other examples of $f$, $\sigma$ and $h$ ensuring 
that $q_t$ evolves in a finite--dimensional family are given in
\cite{benes81a}, \cite{daum86a}. Notice, however, that in these examples
the drift $f$ is not allowed to depend on the observation process $Y$. 

In the present paper we focus on exponential families, according to the 
following
\begin{definition} \label{defiesp} {\bf (Unnormalized exponential family)}
   Let $\{c_1,\cdots,c_m\}$ be scalar functions defined on ${\bf R}$, 
   such that $\{1,c_1,\cdots,c_m\}$ are {\em linearly independent}, 
   have at most
   polynomial growth and are twice continuously differentiable.
   Assume that the convex set 
\begin{displaymath} 
   \Theta_0 := \{\theta\in {\bf R}^m\,:\, 
   \psi(\theta) = \log\; \int \exp[ \theta^T c(x) ]\, d x 
   < \infty \}\ , 
\end{displaymath} 
   has {\em non--empty interior}. 
   Then 
\begin{eqnarray*} 
   EU(c) &=& \{ q(\cdot; \theta,\beta)\,,\, (\theta,\beta) \in U \},  \\ \\
   q(x;\theta,\beta)) &:=& \exp[\theta^T c(x)  + \beta]\ , \
\end{eqnarray*} 
   where  $(\theta,\beta) := 
      [\theta_1, \ \ldots , \ \theta_m, \beta \ ]^T $ and
   $U \subseteq \Theta_0 \times {\bf R}$ is open, 
   is called an unnormalized exponential family of probability densities. 
   The $m+1$ quantities $(\theta,\beta)$ are called the canonical parameters
   for the unnormalized exponential family $EU(c)$. 
\end{definition} 

\begin{remark} \label{remc}
   Given linearly independent   $\{c_1,\cdots,c_m\}$,
   it may happen that 
   the densities $\exp[ \theta^T c(x) ]$ are not integrable. 
   However, it is always possible to extend 
   the family so as to deal with integrable densities only. 
   Indeed, since
   there exist $K > 0$ and $r \geq 0$ such that 
\begin{displaymath} 
       \vert c(x) \vert \leq K\, (1 + \vert x \vert^r)\ , 
\end{displaymath} 
   for all $x\in {\bf R}$, we can define
     $d(x) := \vert x \vert^s$ for all $x\in {\bf R}$, 
   and some $s > r$. Then 
\begin{eqnarray*} 
   EU([c \ \ d])&:=& \{ q'(\cdot; \theta,\mu,\beta)\,,\,\theta\in 
   {\bf R}^m,\mu > 0 \ , \beta \in {\bf R} \}, \\ \\
   q'(x;\theta,\mu, \beta) &:=& \exp[ \theta^T c(x) - \mu\, d(x) 
   + \beta]\ , 
\end{eqnarray*} 
   is an unnormalized exponential family of densities, 
   with a non--empty open parameter set. 
\end{remark} 
In this paper we solve the following problem.
Consider the nonlinear filtering problem for the system
(\ref{Lanc1-1}). Suppose that the nonlinear coefficients
$\sigma_t(\cdot)$ and $h(\cdot)$  are given a priori.
Let be given an exponential family.
Find a drift $f=u$ such that the resulting filtering problem 
has a solution in the given exponential family augmented by
the functions $h$ and $h^2$ in the exponent. More precisely:  
\begin{problem} \label{prob}
Let be given any Lipschitz continuous (uniformly in time) 
diffusion coefficient $a_t(\cdot) = \sigma_t(\cdot)^2$  and any 
Lipschitz--continuous observation function $h$
(which has at most polynomial growth).
% and such that (A) is 
%satisfied for any Lipschitz--continuous (uniformly in time) drift $f$.
%
%{\em (Aggiusta le ipotesi)}
%
Let be given any exponential family $EU(c)$, such that
$EU(c^\bullet) := EU([h \ \ h^2 \ \ c^T]^T)$ is still an exponential family 
according to Definition \ref{defiesp}.
Let $(\zeta,\beta)$ be the  $m +3$
canonical parameters of $EU(c^\bullet)$. 
Set (as comes natural) 
$[\theta_1,\ldots,\theta_m] = [\zeta_3,\ldots,\zeta_{m+2}]$.
Let be given an initial condition $X_0$ such that $\pi_0$ admits an
unnormalized density $q_0 = q(\cdot; \zeta_0, \beta_0)$ 
in the extended family $EU(c^\bullet)$. 
Find a drift $f = u$ such that the filtering problem
for the nonlinear system (\ref{Lanc1-1}) admits a finite dimensional
solution $q_t$ evolving in $EU([h \ \ h^2 \ \  c^T\ ]^T)$.
\end{problem} 
\section{Solution}
In this section we solve Problem \ref{prob} by mean of the following
\begin{theorem} \label{mainres}
 {\bf ( Solution of Problem \ref{prob})}
Assumptions of Problem \ref{prob} in force. 
If the drift
\begin{eqnarray*}
  u_t(x,\zeta_t) &:=& \half \frac{\partial a_t}{\partial x}(x) + 
          \half a_t(x) \zeta_t^T  \frac{\partial c^\bullet}{\partial x}(x),
%          + \\ \nonumber \\ \nonumber
%&& - E_{\theta}\{{\cal L}_t c^\bullet \}^T g^{-1}(\theta) 
%                 \int_{-\infty}^x (c^\bullet (y) - E_{\theta}c^\bullet ) 
%                 \ \exp[\theta^T (c^\bullet (y) - c^\bullet (x))] dy, 
\end{eqnarray*}
with
\begin{equation} 
\zeta^1_t := Y_t + \zeta^1_0, \ \ \zeta^2_t := \zeta^2_0 - \half t,
\ \zeta^i_t := \zeta^i_0, \ i=3,\ldots,m,
\end{equation}
satisfies $(A)$ together with $a_t$ and $h$, then it 
solves Problem \ref{prob}. As a consequence, for the resulting nonlinear
filtering problem with coefficients $u_t, a_t$ and $h$, the optimal
filter is  expressed by the density 
$q_t = q(\cdot,\zeta_t,\beta_0)$  which evolves in the finite dimensional
exponential family $EU(c^\bullet)$.
\end{theorem}
\begin{remark}{\bf (Sufficient conditions for (A))}
For explicit conditions on $u_t, a_t$, $h$ and $X_0$ under which 
(A) holds (and for general results on the solution of the DMZ--SPDE)
see \cite{gyongy90a}. In our case a set of sufficient condition 
ensuring (A) is, as from Theorem 4.1 of \cite{gyongy90a}:
\begin{itemize}
\item[(i)] $h$ is Lipschitz continuous;
\item [(ii)] $a_t, \partial_x a_t, a_t \ \partial_x c^\bullet_j$
are Lipschitz continuous, uniformly in $t$, for $j=2,\ldots,m+2$.
\item[(iii)] $a_t(x) \ y \ \partial_x h(x)$ is Lipschitz continuos
in $(x,y)$, uniformly in $t$. 
\item[(iv)] $Y_0 = 0$ and $q_{X_0} \in EU(c^\bullet)$. 
\end{itemize}
\end{remark}
PROOF of the theorem:
Consider a candidate drift $u_t(\cdot\ ;\zeta_t,\beta_t)$ and the 
associated forward differential operator
\begin{displaymath} 
   {\cal U}_t^\ast \phi = - 
   \frac{\partial}{\partial x}\, [u_t(\cdot\ ;\zeta_t,\beta_t) \ \phi] 
   + \half 
   \frac{\partial^2}{\partial x^2}\, 
   [a_t\, \phi]\ .
\end{displaymath} 
The corresponding state equation originates the filtering problem 
for the system
\begin{eqnarray*}
&&d X_t = u_t(X_t;\zeta_t,\beta_t) dt \ + \sigma_t(X_t) \ \ dW_t, \ \ X_0 \\
&&d Y_t = h(X_t) dt + d V_t\ .
\end{eqnarray*}
%and observations $Y$ as in (\ref{Lanc1-1}).
In order to check under which conditions this filtering problem 
has solution $q_t = q(\cdot\ ;\zeta_t,\beta_t)$ we proceed as follows.
Write the right--hand side of the DMZ--SPDE for the density
$q(\cdot\ ;\zeta_t,\beta_t)$ and equate it to the differential (in time) 
of  $q(\cdot\ ;\zeta_t,\beta_t)$ computed via the chain rule:
\begin{eqnarray*}
  \sum_{i=1}^{m+2}
   \frac{\partial q(\cdot\ ;\zeta_t,\beta_t)}{\partial \zeta_i} 
   \; \circ d {\zeta_t}^i\ 
  +  
   \frac{\partial q(\cdot\ ;\zeta_t,\beta_t)}{\partial \beta} 
   \; \circ d {\beta_t}^i\ 
  \  \\
\  = {\cal U}_t^\ast\, q(\cdot\ ;\zeta_t,\beta_t)\,dt 
   - \half\, q(\cdot\ ;\zeta_t,\beta_t)\, \vert h \vert^2   \,dt 
   +   q(\cdot\ ;\zeta_t,\beta_t)\, h \ \circ dY_t\ .   
\end{eqnarray*}
By dividing both sides by $q$ and straightforward calculations,
one obtains
\begin{eqnarray*}
&& \sum_{i=1}^{m+2} c^\bullet_i  \circ d \zeta^i_t + d \beta_t
=  h  \circ dY_t - \half h^2  dt 
  - u_t(\cdot\ ;\zeta_t,\beta_t) \zeta_t^T
   \frac{\partial c^\bullet}{\partial x}
 \\ && - \frac{\partial u_t(\cdot\ ;\zeta_t,\beta_t)}{\partial x} +
  \frac{1}{2 \ q(\cdot\ ;\zeta_t,\beta_t)} \frac{ \partial^2 [a_t  
      \  q(\cdot\ ;\zeta_t,\beta_t)]}{\partial x^2}.
\end{eqnarray*}
A first step in finding  curves $t \mapsto \zeta_t$ and 
$t \mapsto \beta_t$ such that this equation is satisfied is  to
set (remember that $c^\bullet_1 = h, \ c^\bullet_2 = h^2$)
$d \zeta_t^1 := dY_t, \ \ d \zeta^2_t := - \half dt$.
The above equation reduces then to the (linear) differential equation
\begin{eqnarray*}
\frac{\partial u_t(\cdot\ ;\zeta_t,\beta_t)}{\partial x}
    + \zeta_t^T \frac{\partial c^\bullet}{\partial x}
      \  u_t(\cdot\ ;\zeta_t,\beta_t)  
 \\ \\ =  - \sum_{i=3}^{m+2} \dot{\zeta^i_t} c^\bullet_i - \dot{\beta_t}
 + \frac{1}{2 \ q(\cdot\ ;\zeta_t,\beta_t)} \frac{ \partial^2 [a_t  
      \  q(\cdot\ ;\zeta_t,\beta_t)]}{\partial x^2}.
\end{eqnarray*}
By solving this last equation in $u$ one obtains ($\theta_t$ as
defined in Problem \ref{prob})
\begin{eqnarray*}
u_t(x;\zeta_t,\beta_t) &=& 
\half \frac{\partial a_t(x)}{\partial x} + \half a_t(x)
        \zeta_t^T \frac{\partial c^\bullet}{\partial x} \\ \\
    &-& \dot{\theta_t}^T   
      \exp[- \zeta_t^T c^\bullet(x)] \ \int_{-\infty}^{x}
      c(z) \exp[ \zeta_t^T c^\bullet(z)] dz  \\ \\ 
   && - \dot{\beta_t} \int_{-\infty}^{x}
      \exp[ \zeta_t^T c^\bullet(z)] dz \ .
\end{eqnarray*}
We are free to choose the curve $t \mapsto (\theta_t,\beta_t)$ 
as we wish, as long as it is regular. 
Set $\theta_t := \theta_0$ for all $ t \in [0 \ T] \ $,
and $\beta_t := \beta_0$ for all $t \in [0 \ T] \ $.
\qed
\section{Examples}
We present the following examples of applications
of Theorem \ref{mainres}. 
\begin{example} {\bf Cubic observations.}
The present example is inspired by the cubic sensor 
studied in \cite{HaMaSu}, where it is proven that for the cubic sensor
problem there exists no finite dimensional filter.
Here we consider the case where 
the diffusion coefficient for the state equation is constant
(say $a_t(x) = 1$ for all $t \in [0,\ T]$ and $x \in {\bf R}$)
and with cubic observation function ($h(x) = x^3$ for all
$x \in {\bf R}$). A straightforward application of Theorem
\ref{mainres} yields the following result.

The filtering problem for the system
\begin{eqnarray*}
dX_t &=& [\frac{3}{2} Y_t X_t^2  - 3 (1 + \half t) X_t^5] dt+ \ dW_t,
\ \ q_{X_0} \sim \exp[- x^6], \\ \\
dY_t &=& X_t^3 dt + dV_t, \ \ Y_0=0,
\end{eqnarray*}
is finite dimensional and has conditional law with density
\begin{displaymath}
p_{X_t|{\cal Y}_t}(x) \propto \exp[Y_t x^3 - (1 + \half t) x^6],
\ \ t \in [0, \ T], \ \ x \in {\bf R}.
\end{displaymath}
\end{example}

\begin{example} 
{\bf Linear Observations.} We consider the case of linear observations.
For simplicity, take $h$ equal to the identity function.
Theorem \ref{mainres} yields the following result:
The optimal filter for the filtering problem
\begin{eqnarray*}
dX_t &=& \{\half 
\frac{\partial a_t}{\partial x}(X_t) 
+ \half a_t(X_t) [Y_t + \frac{\mu_0}{v_0} -( t + \frac{1}{v_0})X_t] \}
  dt \\ 
   &+& \sigma_t(X_t) d W_t, \ \ X_0 \sim {\cal N}(\mu_0,v_0), \ \
a_t = \sigma_t^2, \\ \\ 
dY_t &=& X_t dt + dV_t, \ \ Y_0=0,
\end{eqnarray*}
is finite dimensional and has conditional law with density
\begin{displaymath}
X_t|{\cal Y}_t \sim 
{\cal N}(\frac{\mu_0 + Y_t v_0}{1 + v_0 t},\frac{v_0}{1+v_0 t}).
\end{displaymath}
\end{example}

\section{Conclusion}
It seems, at a first sight, that our result contradicts 
classical results on nonexistence of finite dimensional filters,
such as for example Chaleyat--Maurel and Michel (1984) \cite{ChalMaur},
and the related works 
Ocone and Pardoux (1989) \cite{OconPard}, L\'evine (1991) \cite{Levi}.
This contradiction appears a natural consequence of the arbitrariness
of $\sigma$ and $h$. 
Nonetheless, there is no real contradiction.
Indeed, since $\{ \zeta_t, \ \ t \in [0 \ T] \ \}$ depends 
on the observation process $Y$, the drift itself depends on the 
observations. 
This assumption is not allowed in the works mentioned before, and
indeed we cannot construct a nonlinear filtering problem with 
prescribed (nonlinear) $\sigma$ and $h$, with drift $u$
{\em which does not depend on the observation process Y}
and whose solution remains finite dimensional.
We have to allow for observations-dependent drifts in order
to prove our result.
\section{Acknowledgements}
The author did benefit a lot from discussions with Jan H. van Schuppen 
while visiting the CWI institute in Amsterdam.

\end{document}